\documentclass[11pt,a4paper]{article}

\usepackage[utf8]{inputenc} 
\usepackage{geometry} 
\usepackage{subfig} 

\usepackage[english]{babel}
\usepackage{amsmath}
\usepackage{amssymb}
\usepackage{marvosym}
\usepackage{amsfonts}
\usepackage{array}
\usepackage[dvips]{graphicx}

\usepackage{epsfig}
\usepackage{float}
\usepackage{theorem}

\usepackage{amsmath,amssymb,amsfonts,enumerate,rotate,multirow}
\usepackage{bm,latexsym,mathrsfs}
\usepackage{epsfig,graphicx}
\usepackage[usenames]{color}

\newtheorem{thm}{Theorem}[section]
\newtheorem{rem}[thm]{Remark}

\numberwithin{equation}{section}

 \newcommand{\ind}{\mathbf{1}}

 \renewcommand{\tilde}{\widetilde}

\begin{document}

\title{On linear regression models in infinite dimensional spaces with scalar response}
\author{Andrea Ghiglietti$^{a}$, Francesca Ieva$^{a}$, Anna Maria Paganoni$^{b}$, Giacomo Aletti$^{a}$}

\maketitle

\begin{center}
{\small $^a$
       ADAMSS Center \& Department of Mathematics ``F. Enriques'', \\
        Universit\`{a} degli Studi di Milano,\\
        via Saldini 50, 20133, Milan, Italy\\
{\tt andrea.ghiglietti@unimi.it}\\
{\tt francesca.ieva@unimi.it}\\
{\tt giacomo.aletti@unimi.it}\\
\small $^b$
        MOX - Department of Mathematics, \\
        Politecnico di Milano, \\
        via Bonardi 9, 20133, Milan, Italy\\
{\tt anna.paganoni@polimi.it}}
\end{center}


\abstract{In functional linear regression, the parameters estimation involves solving a non necessarily well-posed problem and it has points of contact with a range of methodologies, including statistical smoothing, deconvolution and projection on finite-dimensional subspaces. We discuss the standard approach based explicitly on functional principal components analysis, nevertheless the choice of the number of basis components remains something subjective and not always properly discussed and justified. In this work we discuss inferential properties of least square estimation in this context with different choices of projection subspaces, as well as we study asymptotic behaviour increasing the dimension of subspaces.}

\noindent
\textbf{Keywords:}{Functional Regression, Functional Principal Component Analysis, Asymptotic properties of statistical inference}


\section{Introduction}   \label{introduzione}

It is more and more common in recent years that applications of regression analysis are concerned with functional data. It is the case, for example, when the explanatory variables are curves (or are a digitized points of a curve) linked to a scalar response variable. This arises, for instance, in chemometrics, where some chemical variable has to be predicted by a digitized signal such as the Near Infrared Reflectance (NIR) spectroscopic information (see \cite{Osborneetal1984}, \cite{FrankFriedman1993}). Other examples concern environmental problems, like prediction of total annual precipitation for Canadian weather stations from the pattern of temperature variation through the year \cite{RS2006}, or linguistic issues \cite{Hastieetal1993}, like the analysis of the relationship between log-spectra of sequences of spoken syllables and phoneme classification \cite{MarxEilers1996}.

In all these cases, classical regression models for multivariate data may be inadequate, since the functional nature of covariates should be exploited using proper estimation and inferential techniques.

In other words, in functional linear regression, the parameters estimation involves solving an illposed problem \cite{Cuevasetal2002} and has points of contact with a range of methodologies, including statistical smoothing and deconvolution (see, among others, \cite{Cardotetal2003} and references therein). The standard approach to carry out estimation and inference on regression parameters is based explicitly on functional principal components analysis (FPCA, see \cite{RS2006} and references therein) and, consequently, on spectral decomposition in terms of eigenvalues and eigenfunctions. Despite FPCA, or analogous projection methods, are often effective and straightforward to apply to the analysis of functional data, the choice of the number of basis components remains something subjective and not always properly discussed and justified. Even if several criteria exist to determine the number of basis functions to be selected dimensional reduction methods per se do not ensure the proper estimation of the regression parameters. We show that, given the sub-space identified by the the chosen basis, the classical procedures do not automatically ensure to obtain an unbiased estimate neither of the true functional coefficient nor of its projection on the correspondent sub-space.

In this work we face the functional linear model with scalar response. In our model a real random response $Y$ is linked to a square integrable random function X defined on some compact set $T$ of $\mathbb{R}$, as

\begin{equation}
y_i=\int_T x_i(t)\beta(t)dt + \epsilon_i,\ \ \ i=1,..,n,
\end{equation}\label{eq:the_model_functional}

\noindent We discuss the choice of suitable finite sub-spaces of $L^2(T)$, called identifiable sub-spaces, where the least square estimation problem is well posed. We point out the properties in terms of bias and variance of the related estimators.
Moreover we explain the reasons why the FPCA comes out to be the optimal solution of a bias-variance trade off problem when no information are available on the space where the regression parameters are defined. Finally we discuss the influence on the parameters estimates (in terms of bias) of the orthogonal component of the sub-space identified by the FPC basis, and we provide a simulation study that shows the theoretical results.

The paper is organized as follows: firstly, the model setting and the functional parameters estimation (Section \ref{section_model}) together with a critical discussion of their inferential properties (Section \ref{section_estimation_subspace}) are presented for finite dimensional sub-spaces. Then the large dimensional case is considered (Section \ref{section_D_inf}) and asymptotic results for increasing size of the sub-space are introduced.
Appendixes \ref{appendix_relation_D_E} and \ref{appendix_cauchy} gather some auxiliary results while Appendix \ref{appendix_simulation} contains the setting details of the simulation study. All the analyses are carried out with \texttt{R}, see \cite{R}.

\section{Model setting and functional least square estimation}   \label{section_model}

Let consider the functional model in \eqref{eq:the_model_functional},
\begin{equation*}
y_i=\int_T x_i(t)\beta(t)dt + \epsilon_i,\ \ \ i=1,..,n,
\end{equation*}
where $\beta\in L^2(T)$, with $T$ compact set of $\mathbb{R}$, $x_i\in L^2(T)$ and $\epsilon_i\in\mathbb{R}$.
We will consider $\beta$ as a deterministic function and $\epsilon_1,..,\epsilon_n$ as random variables independent of $x_1,..,x_n$, with $\bm{E}\left[\epsilon_i\right]=0$ and $\bm{Var}\left[\epsilon_i\right]=\sigma^2>0$.
We assume to collect the values of the outcomes $y_1,..,y_n$ and to observe the data $x_i$ only in correspondence of $p$ discrete values $t_1,..,t_p\in T$,
i.e. the data are $\left(x_i(t_1),..,x_i(t_p),y_i\right)$ for $i=1,..,n$.
This is the case treated for example in \cite{Osborneetal1984,FrankFriedman1993,Hastieetal1993,MarxEilers1996,RS2006}, among others.
To ease notation, we let $\langle a,b\rangle$ denote the usual inner product in $L^2(T)$, as $\langle a,b\rangle=\int_Ta(t)b(t)dt$,
and $\|a\|$ is the corresponding norm $\sqrt{\int_Ta^2(t)dt}$.
Accordingly, we can write the model~\eqref{eq:the_model_functional} as
\begin{equation}\label{eq:the_model}
y_i=\langle x_i,\beta\rangle + \epsilon_i,\ \ \ i=1,..,n.
\end{equation}
To obtain the asymptotic results presented in the paper,
we assume that $x_1,..,x_n$ are i.i.d. realizations of a process $X$ with support $S_X\subseteq L^2(T)$,
zero mean and bounded second moment,
i.e. $\bm{E}[X]=0$ and $\bm{E}[\|X\|^2]<\infty$.
In general, neither the distribution of the random process $X$ nor its support $S_X$ are assumed to be known.
The quantities $\epsilon_1,..,\epsilon_n$ model the errors in observing the outcomes $y_1,..,y_n$, and so are assumed to be unknown.
The function $\beta$ is unknown and its estimation is the main focus of this paper.

We need the following setting to describe the functional estimation presented in the paper:
let $S$ be the smallest closed sub-space of $L^2(T)$ such that $S_X\subseteq S$,
and we call $S^{\perp}$ the sub-space of $L^2(T)$ orthogonal to $S$, so that
\begin{equation}\label{eq:decomposition_L2_1}
L^2(T)=S\oplus S^{\perp}.
\end{equation}
In general, the set $S$ may not coincide with $L^2(T)$, that means $S^{\perp}\neq \emptyset$.
For instance, consider the following process in $L^2(T)$, with $T=[-1,1]$:
\begin{equation}\label{eq:example_X}
X(t)\ =\ \sum_{k=0}^{\infty} U_k \eta_k \varphi_k(t),
\end{equation}
where $\{U_k:k\geq0\}$ are i.i.d. uniform random variables in $[-1,1]$,
$\{\eta_k:k\geq0\}$ a sequence of positive coefficients such that $\sum_k\eta_k^2<\infty$,
$\varphi_0(t) = 1/\sqrt{2}$ and $ \varphi_k(t)=\cos{(\pi k t)}; k\geq1.$

In this case, the support $S_X$ is composed by the even functions such that
$|\langle g,\varphi_k\rangle|\leq \eta_k$, for any $k\geq0$ and $g\in L^2(T)$.
Then, the smallest sub-space of $L^2(T)$ including $S_X$ coincides with the set of the even functions,
while the orthogonal space is represented by the odd functions, i.e.
$$S:=\left\{\varphi_k(t);k\geq0\right\},\ \ \ \ \ \ \ S^{\perp}:=\left\{\sin{(\pi k t)};k\geq1\right\}.$$

\begin{rem}\label{rem:intercetta}
It is worth saying that the results on the estimation of $\beta$ presented in the paper also hold when the model is
$$y_i=\alpha+\langle x_i,\beta\rangle + \epsilon_i,\ \ \ i=1,..,n,$$
with $\alpha\in\mathbb{R}$, or when $\bm{E}[X]\neq0$.
In these cases, the model~\eqref{eq:the_model} is applied to the centered data, i.e. $y_i-(\sum_{i=j}^ny_i)/n$ and $x_i-(\sum_{i=j}^nx_i)/n$
$i=1,..,n$, so that the asymptotic results are straightforwardly verified.
\end{rem}

\subsection{Functional least square estimation in finite sub-spaces}   \label{section_least square}

In the multivariate regression analysis, a common approach to solve the problem of the estimation of $\beta$ is to compute the least square estimator.
However, it is well known that this approach can't be straightforwardly generalized to the functional context,
not even in the case of $x_1,..,x_n$ entirely observed for any $t\in T$.
In fact, the extension of the least square estimator to the functional framework would be
\begin{equation}\label{eq:minimum_problem_complete}
\widehat{\beta}_n:=\arg\min_{b\in L^2(T)}\ f(b)=\arg\min_{b\in L^2(T)}\ \left\{\sum_{i=1}^n(y_i-\langle x_i,b\rangle )^2\right\}
\end{equation}
and it is trivial to note that for any $n\in\mathbb{N}$, $x_1,..,x_n\in L^2(T)$ and $y_1,..,y_n\in \mathbb{R}$ there exist infinite functions $b\in L^2(T)$ such that $f(b)=0$, even for $S=L^2(T)$. Then, the estimator $\widehat{\beta}_n$ can never be well defined by following the least square approach~\eqref{eq:minimum_problem_complete}.

However, a least square estimator of $\beta$ can be computed in a finite sub-space $D$ of $ L^2(T)$.
In fact, let $D\subset L^2(T)$ be a sub-space where the data $x_i$ are reconstructed from their discrete observation $\left(x_i(t_1),..,x_i(t_p)\right)$ by classical smoothing techniques, so obtaining $x_i^D\in D$.
Therefore, we will simply assume that $x_i^D$ represents the projection of $x_i$ on $D$, and in particular that
\begin{equation}\label{eq:reconstruction_assumption}
\langle x_i^D,g\rangle\ =\ \langle x_i,g\rangle,\ \ \ \forall g\in D.
\end{equation}
Given~\eqref{eq:reconstruction_assumption}, the following minimization problem is, under mild conditions, well posed:
\begin{equation}\label{eq:minimum_problem_D}
\widehat{\beta}^D_n:=\arg\min_{b\in D}\ f(b)=\arg\min_{b\in D}\ \left\{\sum_{i=1}^n(y_i-\langle x_i,b\rangle )^2\right\},
\end{equation}
and it can be computed exactly since from~\eqref{eq:reconstruction_assumption} the real function $x_i\in S$ can be replaced in~\eqref{eq:minimum_problem_D} by its reconstruction $x_i^D\in D$.
First, note that if $d:=\dim(D)$ is greater than the sample size $n$, or than the number $p$ of observation points,
the solution of~\eqref{eq:minimum_problem_D}
is not unique as in~\eqref{eq:minimum_problem_complete}, which provides us the condition $d\leq \min\{n;p\}$.
Moreover, if there exists $\beta_0\in D\cap S^{\perp}$ then $\langle x_i,\beta^D+\beta_0\rangle=\langle x_i,\beta^D\rangle$ for any $\beta^D\in D$, which implies that the minimum is not unique and so $\widehat{\beta}^D_n$ is not well defined.
From~\eqref{eq:reconstruction_assumption} we have that the same situation occurs when we replace $x_i$ with its reconstruction $x_i^D$.
To avoid this problem, we introduce the following concept.

We call \emph{identifiable} any sub-space $D$ such that $D\cap S^{\perp}=0$.
We recall that $S$ and $S^{\perp}$ are individuated by $S_X$, which is in general unknown.
Then, the statistician has the important role of choosing a sub-space with no components orthogonal with respect to the sample data,
which are formally the components lying in $S^{\perp}$.

It is worth highlighting that estimating $\beta$ in a finite sub-space $D$
is intrinsically a consequence of the reconstruction procedure of the data $x_i$ on $D$.
In fact, if we consider the problem~\eqref{eq:minimum_problem_complete} computed with the reconstructed data $x_i^D\in D$,
it is easy to see that for any $b_1,b_2$ such that $(b_1-b_2)\in D^{\perp}$,
we have $f(b_1)=f(b_2)$ and the solution of~\eqref{eq:minimum_problem_D} can't be unique.
Hence, the uniqueness of the solution of~\eqref{eq:minimum_problem_D}
can be obtained only by restricting the problem to the sub-space where the data have been reconstructed, that is $D$ is an identifiable sub-space.

Moreover, in some application, the physical context of the problem may provide a prior information on $\beta$ so that
searching a solution in a specific sub-space $D$ could be the smartest thing to do.
In this case, even if the data $x_i$ are perfectly recorded at any $t\in T$, the problem~\eqref{eq:minimum_problem_D}
would only consider their projection on $D$, since the part on $D^{\perp}$ is useless.
In fact, from~\eqref{eq:minimum_problem_D} the components of the data $x_i$ orthogonal to $D$ are irrelevant,
because $b^D\in D$ and the orthogonal part vanishes in the scalar product $\langle x_i,b^D\rangle $.
Then, the strategy of searching $\beta\in D$ through~\eqref{eq:minimum_problem_D} suggests to reconstruct the data on $D$.

In practice, the a priori information on $\beta$ may not guarantee to determine a finite sub-space $D$ where $\beta$ belongs to.
Then, the sub-space $D$ is typically chosen to reconstruct the data $x_1,..,x_n$ at best,
and so we can imagine that in general the true $\beta$ may not lie in that sub-space $D$.
In this case, it is not clear what $\widehat{\beta}^D_n$ defined in~\eqref{eq:minimum_problem_D} is actually estimating,
and which are its statistical properties.
In the following section, we provide an answer to this issue.
For instance, we will show that, in general, the least square estimator $\widehat{\beta}^D_n$ computed on $D$
does not converge to the projection of the real $\beta$ on $D$, as one may expect.
Moreover we will discuss the collinearity effects in the estimation of $\beta$,
which plays a central role in the unbiasedness and consistency of the estimator $\widehat{\beta}^D_n$.

\section{Properties of Least Square Estimator in finite sub-spaces}   \label{section_estimation_subspace}

To investigate the statistical properties of $\widehat{\beta}^D_n$,
we rewrite~\eqref{eq:minimum_problem_D} in a slightly different way.
First, we introduce the projection operator $\pi:D\rightarrow S$ of $D$ on $S$, and call $E\subset S$ the image of $\pi$, i.e.
$$E\ :=\ \left\{\ x\in S\ :\ \exists y\in D, x=\sum_{k=1}^{\dim(S)}\langle y,\varphi^S_k\rangle\varphi^S_k\ \right\},$$
where $\left\{\varphi_k^S;k=1,..,\dim(S)\right\}$ is an orthonormal basis of $S$.
Naturally, the definition of $E$ implies that $D\subset E \oplus S^{\perp}$.

Appendix~\ref{appendix_relation_D_E} is dedicated to explore more precisely the relation among $D$ and $E$:
for any given $D$ and $S$, we describe how to compute an orthonormal basis for $E$ and we provide the analytic expression of the projection operator $\pi$.
Here, we focus on the following properties:
\begin{itemize}
\item[(1)] since $D\cap S^{\perp}=0$ ($D$ is identifiable) and $D$ is finite dimensional, it is possible to show that $\pi$ is invertible (see Appendix~\ref{appendix_relation_D_E}), so that $\pi$ is a bicontinuous operator from $D$ to $E$;
\item[(2)] for any $b^D\in D$, calling $b^{E}=\pi(b^D)$, we have that
$$\langle x_i,b^D\rangle=\langle x_i,b^{E}\rangle+\langle x_i,b^D-b^{E}\rangle=\langle x_i,b^{E}\rangle,$$
because $(b^D-b^{E})\in S^{\perp}$ and $x_i\in S$.
\end{itemize}
From (1) and (2), we have that $f(b^D)=f(b^{E})$ for any $b^D\in D$, so that the element of $D$ that minimizes $f$ is univocally associated
through the projection $\pi$ with the element of $E$ minimizing $f$.
Hence, the least square estimator computed minimizing in~\eqref{eq:minimum_problem_D} can be obtained as
\begin{equation}\label{eq:inverse_projection}
\widehat{\beta}^{D}_n=(\pi)^{-1}(\widehat{\beta}^{E}_n),
\end{equation}
where
\begin{equation}\label{eq:minimum_problem_E}
\widehat{\beta}^{E}_n:=\arg\min_{b\in E}\ f(b)=\arg\min_{b\in E}\ \left\{\sum_{i=1}^n(y_i-\langle x_i,b\rangle )^2\right\}.
\end{equation}
Then, in the following, we study the statistical properties of $\widehat{\beta}^{E}_n$
to describe the behavior of the estimator $\widehat{\beta}^{D}_n$ computed in $D$,
which is the sub-space individuated by the experimenter as mentioned before.
The problem~\eqref{eq:minimum_problem_E} is solved in Subsection~\ref{section_estimation_betaE},
where the properties of $\widehat{\beta}^{E}_n$ are investigated.
After that, a wide analysis on the behavior $\widehat{\beta}^{D}_n$ is detailed in Subsection~\ref{section_a_b_c}.

Finally, to sake of simplicity, we define the sub-space $F=S\cap E^{\perp}$, so that we replace~\eqref{eq:decomposition_L2_1}
with the following expression
\begin{equation}\label{eq:decomposition_L2_2}
L^2(T)=E\oplus F\oplus S^{\perp}.
\end{equation}
Then, a unique orthogonal decomposition can be realized for any $\beta\in L^2(T)$:
\begin{equation}\label{eq:decomposition_beta}
\beta\ =\ \beta^E\ +\ \beta^F\ +\ \beta^{S^{\perp}},
\end{equation}
where $\beta\in D$ implies $\beta^F=0$, since $D\subset E \oplus S^{\perp}$.

\subsection{Characterization of the least square estimator in $E$}   \label{section_estimation_betaE}

In this section, we focus on solving~\eqref{eq:minimum_problem_E} and
we obtain the main properties of $\widehat{\beta}^{E}_n$.
Given any orthonormal basis for $D$ and $S$,
denoted by $\left\{\varphi_k^D;k=1,..,d\right\}$ and $\left\{\varphi_k^S;k=1,..,\dim(S)\right\}$ respectively,
we can compute the orthonormal basis for $E$ and we denote it as $\bm{\varphi^{E}(t)}:=\left\{\varphi_k^E;k=1,..,d\right\}$
(see Appendix~\ref{appendix_relation_D_E} for the details).
Then, we call $x_i^{E}$ the projection of $x_i$ on $E$ and note that for any $b^{E}\in E$
$$\langle x_i,b^{E}\rangle=\langle x_i^{E},b^{E}\rangle+\langle x_i-x_i^{E},b^{E}\rangle=\langle x_i^{E},b^{E}\rangle,$$
since $(x_i-x_i^{E})$ lies in a sub-space orthogonal to $E$.
Hence,~\eqref{eq:minimum_problem_E} can be solved with finite dimensional quantities,
obtaining $\widehat{\beta}^{E}_n(t):=(\bm{\widehat{\beta}^{E}_n})^T
\cdot\bm{\varphi^{E}(t)}$, where
\begin{equation}\label{eq:minimum_problem_E_smart}
\bm{\widehat{\beta}^{E}_n}:=\arg\min_{\bm{b^{E}}\in \mathbb{R}^d}\
\left\{(\bm{y}-X^{E}\bm{b^{E}})^T(\bm{y}-X^{E}\bm{b^{E}})\right\},
\end{equation}
$\bm{y}$ is the $n$-vector composed by the observed values $\bm{y}=(y_1,..,y_n)^T$ and
$X^{E}$ is the $n\times d$-matrix, where $[X^{E}]_{ij}=\langle x_i,\varphi_j^{E}\rangle$.
As in the multivariate theory, we can easily obtain
\begin{equation}\label{eq:estimate_E}
\bm{\widehat{\beta}^{E}_n}=((X^{E})^TX^{E})^{-1}(X^{E})^T\bm{y}.
\end{equation}
Now, let us discuss the statistical properties of $\bm{\widehat{\beta}^{E}_n}$.
Using decomposition~\eqref{eq:decomposition_beta}, the model~\eqref{eq:the_model} can be written as
$$y_i=\langle x_i,\beta\rangle + \epsilon_i=\langle x_i,\beta^{E}\rangle+\langle x_i,\beta^{F}\rangle + \epsilon_i$$
for any $i=1,..,n$, since $x_i$ is orthogonal to $\beta^{S^{\perp}}$.
In matrix notation, the last expression becomes
\begin{equation*}
\bm{y}=\langle \bm{x},\beta^{E}\rangle+\langle \bm{x},\beta^{F}\rangle+ \bm{\epsilon}
\end{equation*}
where $\bm{x(t)}=(x_1(t),..,x_n(t))^T$, $\bm{\epsilon}=(\epsilon_1,..,\epsilon_n)^T$ and $\bm{y}=(y_1,..,y_n)^T$.
Since the dimension of $E$ is finite, we can rewrite the last expression as follows
\begin{equation}\label{eq:espression_collinearity}
\bm{y}=X^{E}\bm{\beta^{E}}+\langle \bm{x},\beta^{F}\rangle+ \bm{\epsilon},
\end{equation}
where $\bm{\beta^{E}}\in\mathbb{R}^d$ is the vector such that $\beta^{E}(t)=(\bm{\beta^{E}})^T\cdot\bm{\varphi^{E}(t)}$.
Note that the estimator $\bm{\widehat{\beta}^{E}_n}$ is computed in~\eqref{eq:estimate_E} only with the data projected on $E$, i.e. $X^{E}$;
then, the quantity $\langle \bm{x},\beta^{F}\rangle$ in~\eqref{eq:espression_collinearity} represents the part of the data
which has not been used to compute $\bm{\widehat{\beta}^{E}_n}$, so that the least square estimation approach treats $\langle \bm{x},\beta^{F}\rangle$ in~\eqref{eq:espression_collinearity} as the independent error $\bm{\epsilon}$.
Nevertheless, the quantity $\langle \bm{x},\beta^{F}\rangle$ can be correlated to $X^{E}$, and this correlation
plays a central role in the estimation of $\beta$.

To characterize $\bm{\widehat{\beta}^{E}_n}$, we substitute~\eqref{eq:espression_collinearity} in~\eqref{eq:estimate_E}, obtaining
\[\begin{aligned}
\bm{\widehat{\beta}^{E}_n}\ &&=&\ \bm{\beta^{E}}+((X^{E})^TX^{E})^{-1}(X^{E})^T\langle \bm{x},\beta^{F}\rangle+
((X^{E})^TX^{E})^{-1}(X^{E})^T\bm{\epsilon}\\
&&=&\ \bm{\beta^{E}}+\bm{\gamma_n^d}+((X^{E})^TX^{E})^{-1}(X^{E})^T\bm{\epsilon},\\
\end{aligned}\]
where $\bm{\gamma_n^d}:=((X^{E})^TX^{E})^{-1}(X^{E})^T\langle \bm{x},\beta^{F}\rangle$.
Then, conditioning to the data $\bm{x(t)}$, the quantity $\bm{\widehat{\beta}^{E}_n}$ presents the following features:
\begin{equation}\label{eq:features_beta_E}
\bm{E}\left[\bm{\widehat{\beta}^{E}_n}|\bm{x}\right]=\bm{\beta^{E}}+\bm{\gamma_n},\ \ \ \ \
\bm{Cov}\left(\bm{\widehat{\beta}^{E}_n}|\bm{x}\ \right)=\sigma^2((X^{E})^TX^{E})^{-1}.
\end{equation}
The term $\bm{\gamma_n}$ catches the relation among $X$ on $E$ and $X$ on $\beta^F$,
see also Subsection~\ref{section_a_b_c}.
Moreover, since $x_1,..,x_n$ are i.i.d. realizations of $X$ and $\bm{E}[\|X\|^2]<\infty$,
we can apply the Strong Law of Large Numbers (SLLN) obtaining
\begin{equation}\label{eq:convergence_bias}
\bm{\gamma_n}\ \stackrel{a.s.}{\rightarrow}\ \bm{\gamma}:=\left(\bm{E}\left[\bm{X^{E}}(\bm{X^{E}})^T\right]\right)^{-1}
\bm{E}[\bm{X^{E}}\langle X,\beta^{F}\rangle],
\end{equation}
where $\bm{X^{E}}:=\langle X,\bm{\varphi^{E}(t)}\rangle\in\mathbb{R}^d$.
Using~\eqref{eq:convergence_bias} we get that $\bm{\widehat{\beta}^{E}_n}\stackrel{a.s.}{\rightarrow}_n\bm{\beta^{E}}+\bm{\gamma}$.
The quantity $\bm{\gamma}$ has a direct functional representation given by $\gamma(t)=(\bm{\gamma})^T\cdot\bm{\varphi^{E}(t)}$,
and we directly obtain the consistency of $\widehat{\beta}^{E}_n$:
\begin{equation}\label{eq:convergence_beta_E}
\widehat{\beta}^{E}_n\ \stackrel{a.s.}{\rightarrow}\ \beta^{E}\ +\ \gamma.
\end{equation}

\begin{rem}\label{rem:bias_delta}
Note that, since $\bm{E}[X]=0$, the bias $\bm{\gamma}$ can also be written as
\begin{equation}\label{eq:bias_zero_mean}
\bm{\gamma}\ =\ \left(\Sigma^E\right)^{-1}
\bm{Cov}(\bm{X^{E}},\langle X,\beta^{F}\rangle),
\end{equation}
where $\Sigma^E:=\bm{Cov}\left[\bm{X^{E}}\right]$.
The meaning of $\bm{\gamma}$ can be easily seen when it is represented along the principal components of $\bm{X^{E}}$.
If we denote with $V^E$ the matrix composed by the eigenvectors ($\bm{\psi_1},..,\bm{\psi_d}$) of $\Sigma^E$,
and if we call $Z^{E}_k:=(\bm{X^{E}})^T\cdot \bm{\psi_k}=\langle X,\psi_k\rangle$ for $k=1,..,d$,
the bias along the $k^{th}$ principal components (i.e. $\delta_k=\langle \gamma,\psi_k\rangle$) can be express as follows
$$\delta_k\ =\ \left(\bm{Var}\left[Z^{E}_k\right]\right)^{-1}\cdot
\bm{Cov}\left[Z^{E}_k,\langle X, \beta^{F}\rangle\right]=
\left(\frac{\bm{Var}\left[\langle X, \beta^{F}\rangle\right]}{\bm{Var}\left[Z^{E}_k\right]}\right)^{1/2}\cdot
\bm{Cor}\left[Z^{E}_k,\langle X, \beta^{F}\rangle\right],$$
which shows how the bias reflects the correlation among $X$ on $E$ and $X$ on $F$.
\end{rem}

\subsection{Discussion on the least square estimation in finite identifiable sub-spaces}   \label{section_a_b_c}

We now discuss the behavior of the least square estimator in the finite sub-space $D$, i.e. $\widehat{\beta}^{D}_n$.
To this aim, we consider the results~\eqref{eq:features_beta_E} and~\eqref{eq:convergence_beta_E} related to $\widehat{\beta}^{E}_n$,
and, through the relation~\eqref{eq:inverse_projection}, we discuss the properties of $\widehat{\beta}^{D}_n$.
In particular, in this subsection we focus on the asymptotic behavior of $\widehat{\beta}^{D}_n$,
even if analogous arguments can be used to describe its bias for fixed $n$.
Since $\pi^{-1}:E\rightarrow D$ is continuous,
the consistency of $\widehat{\beta}^{D}_n$ can be easily obtained from~\eqref{eq:convergence_beta_E}:
\begin{equation}\label{eq:convergence_beta_D}
\widehat{\beta}^{D}_n\ \stackrel{a.s.}{\rightarrow}\ \pi^{-1}\left(\beta^{E}+\gamma\right).
\end{equation}
The real issue here is to understand what this limit represents.
The discussion is structured as follows: we analyze the consistency of the least square estimator $\widehat{\beta}^{D}_n$
in these different cases
\begin{itemize}
\item[(a)] $\beta\in D$;
\item[(b)] $\beta\notin D$, but $\beta^F=0$;
\item[(c)] $\beta\notin D$, and $\beta^F\neq0$.
\end{itemize}

\textit{Case (a): $\beta\in D$.} In this situation, we trivially have $\beta^F=0$ since $D\subset E \oplus S^{\perp}$;
this implies $\beta=\beta^E+\beta^{S^{\perp}}$ and
$\gamma_n=\gamma=0$ for any $n\geq1$ by definition.
Moreover, since $(\beta^E+\beta^{S^{\perp}})\in D$, we have that $\pi^{-1}\left(\beta^{E}\right)=\beta^E+\beta^{S^{\perp}}$.
Hence, we obtain
$$\|\widehat{\beta}^{D}_n-\beta\|\ \stackrel{a.s.}{\rightarrow}\ 0.$$
Then, when the true $\beta$ belongs to the sub-space $D$, the least square estimator on $D$ is consistent.
In Figure~\ref{Simulazione_1}-third panel, we report $100$ independent simulations detailed in Appendix~\ref{appendix_simulation}
in which an estimate of $\widehat{\beta}^{D}_n$ is computed for large $n$ and $\beta \in D$.
The pointwise mean of the estimates of $\widehat{\beta}^{D}_n$ (dotted line) is very close to the true $\beta$ (solid line).
This shows that the estimator $\widehat{\beta}^{D}_n$ is unbiased and consistent.
\begin{figure}
\centering
  \includegraphics[width=0.4\textwidth]{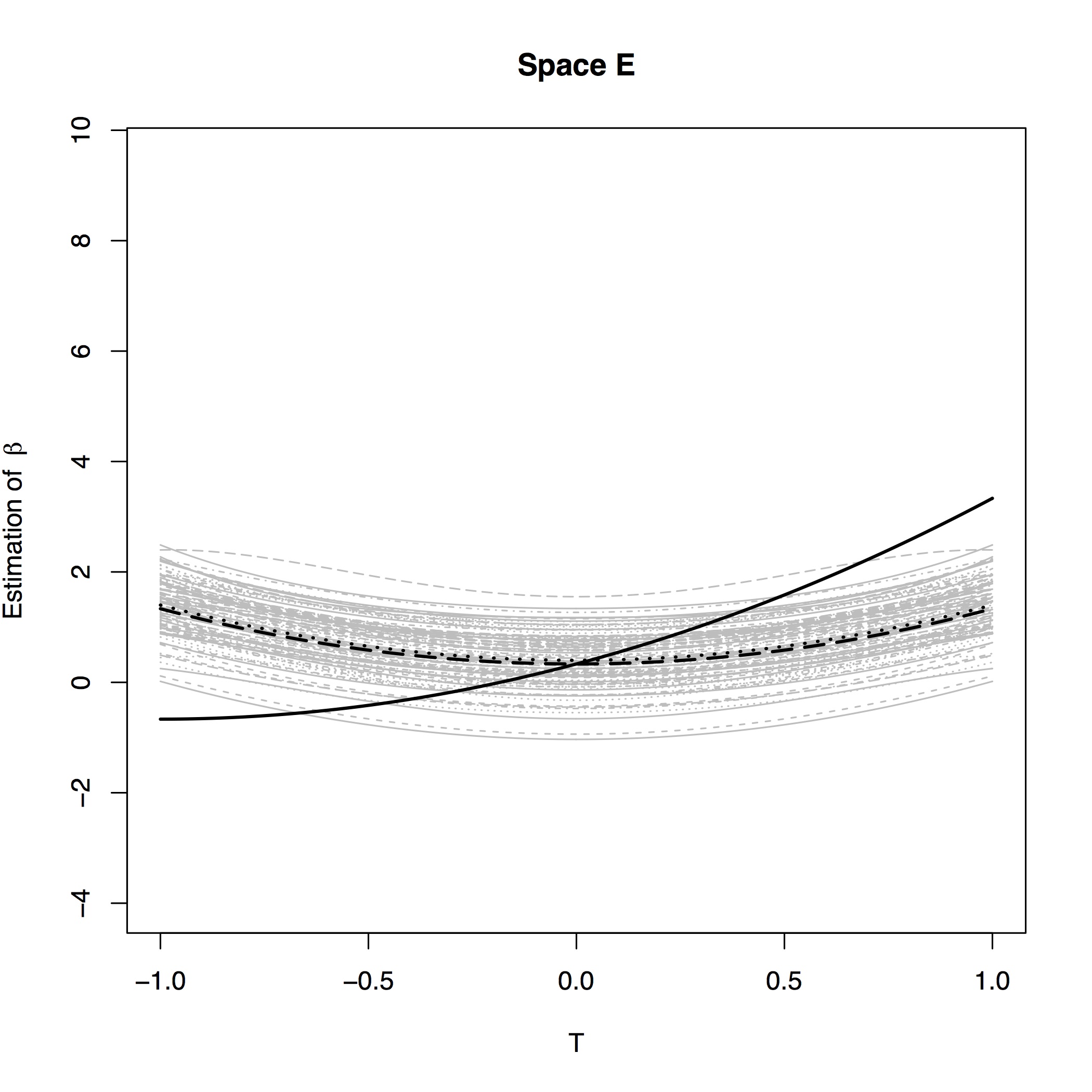}
  \includegraphics[width=0.4\textwidth]{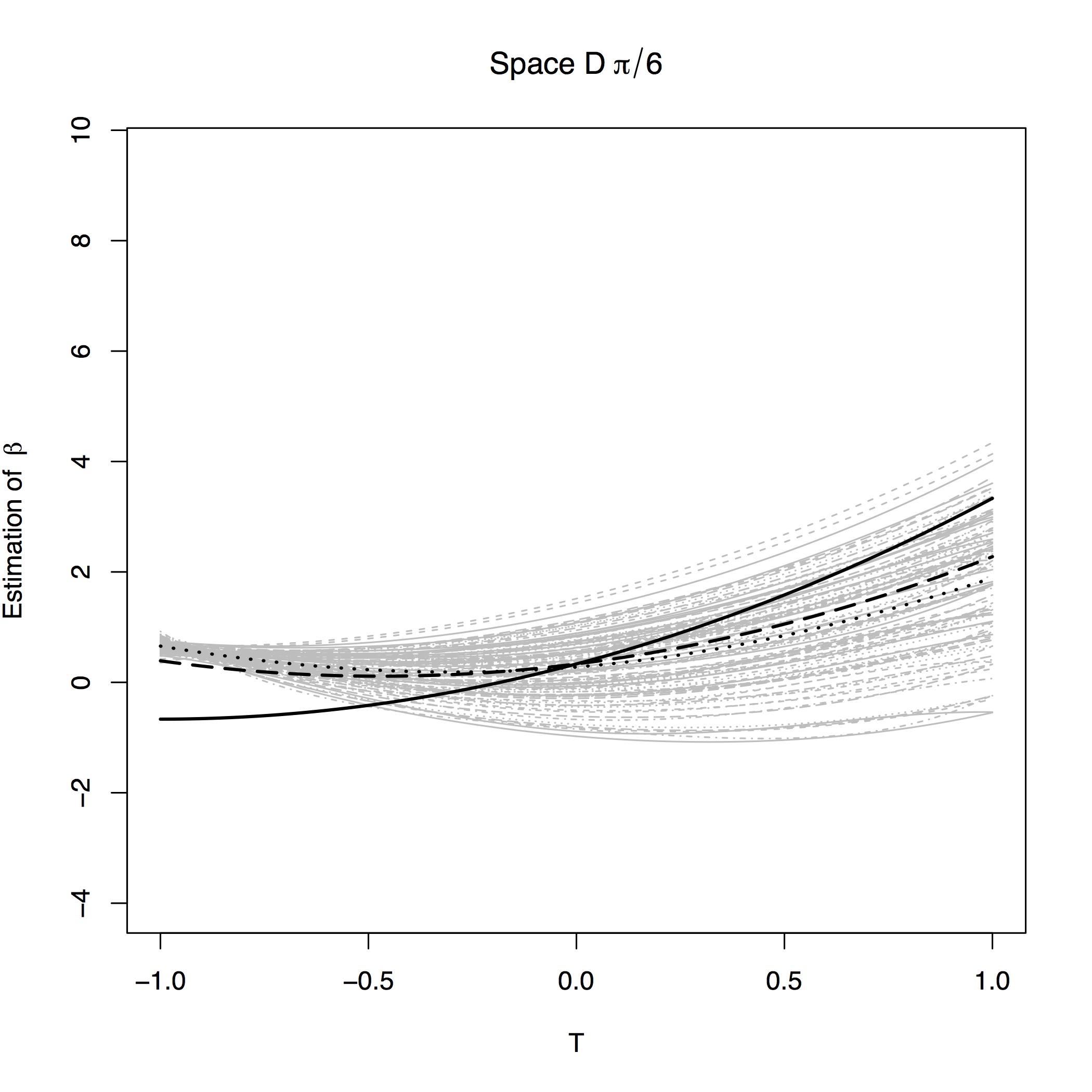}
  \includegraphics[width=0.4\textwidth]{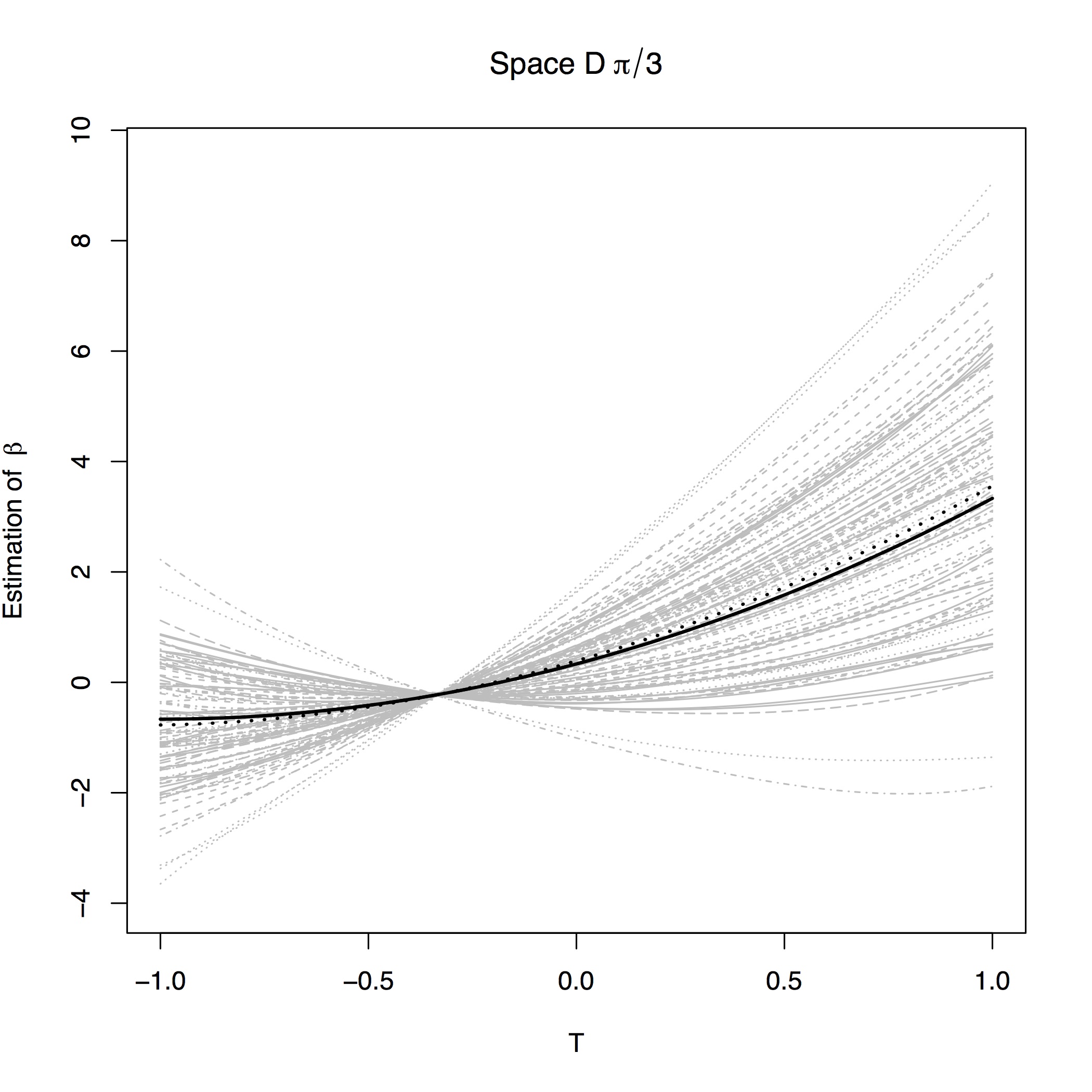}
    \caption{The solid lines are the true $\beta(t)$, the dashed lines are the projection of $\beta(t)$ in the corresponding spaces,
    the dotted lines are the pointwise means of the estimated $\hat\beta(t)$.
    For further details on the simulation setting see Appendix~\ref{appendix_simulation}.}
    \label{Simulazione_1}
\end{figure}\\

\textit{Case (b): $\beta\notin D$, but $\beta^F=0$.}
Analogously to case (a), $\beta^F=0$ implies $\beta=\beta^E+\beta^{S^{\perp}}$ and
$\gamma_n=\gamma=0$ for any $n\geq1$.
However, in this case $\pi^{-1}\left(\beta^{E}\right)\neq\beta^E+\beta^{S^{\perp}}$, so that
$$\|\widehat{\beta}^{D}_n-\beta\|\ \stackrel{a.s.}{\rightarrow}\ \|\left(\pi^{-1}\left(\beta^{E}\right)-\beta^E\right)-\beta^{S^{\perp}}\|,$$
which means the estimator $\widehat{\beta}^D_n$ is not consistent for $\beta$.
The asymptotic bias belongs to the sub-space orthogonal to the data, i.e.
$$\left(\pi^{-1}\left(\beta^{E}\right)-\beta^E\right)-\beta^{S^{\perp}}\ \in\ S^{\perp}.$$
This latter fact can be seen in Figure~\ref{Simulazione_1}-first and second panels.
In fact, in this simulation the difference among $\beta$ (solid line) and the pointwise means of the estimates of $\widehat{\beta}^{D}_n$
(dotted lines) are odd functions, i.e. $S^{\perp}$ in the example.

Since the errors in estimating $\beta$ with $\widehat{\beta}^D_n$ belongs to a space which can't be explored by the data,
the bias can be eliminated only by using a priori information on $\beta$ to modify the choice of $D$.
It is also worth observing that this bias is totally irrelevant if the interest in estimating $\beta$ is only related to the quantity $\langle x,\beta\rangle$ in the regression context,
because the estimation and inference of the inner product is not influenced by any component of $\beta$ in $S^{\perp}$.

In other words in all the cases of Figure~\ref{Simulazione_1} the pointwise means of the estimates of $\widehat{\beta}^{D}_n$ (dotted lines)
only differs in their odd component. Hence, since the data are even functions, the inference on $\langle X,\beta\rangle$ is equivalent.
Summing up in case (b), the choice of $D$ does not influence the explanation of the phenomena related to the regression,
but it is relevant when the interest lies in the reconstruction of the true $\beta$.\\

%

\begin{figure}
\centering
  \includegraphics[width=0.4\textwidth]{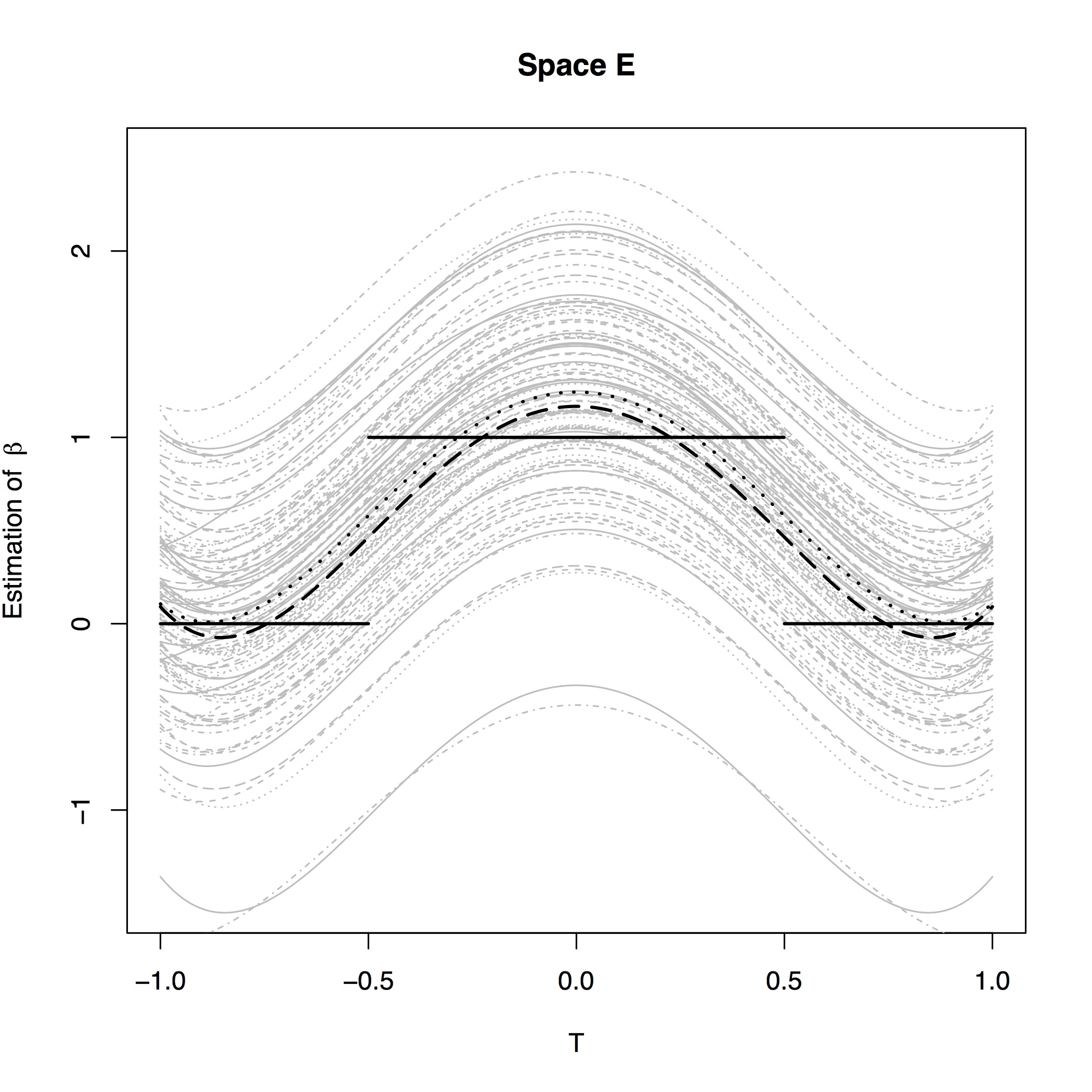}
  \includegraphics[width=0.4\textwidth]{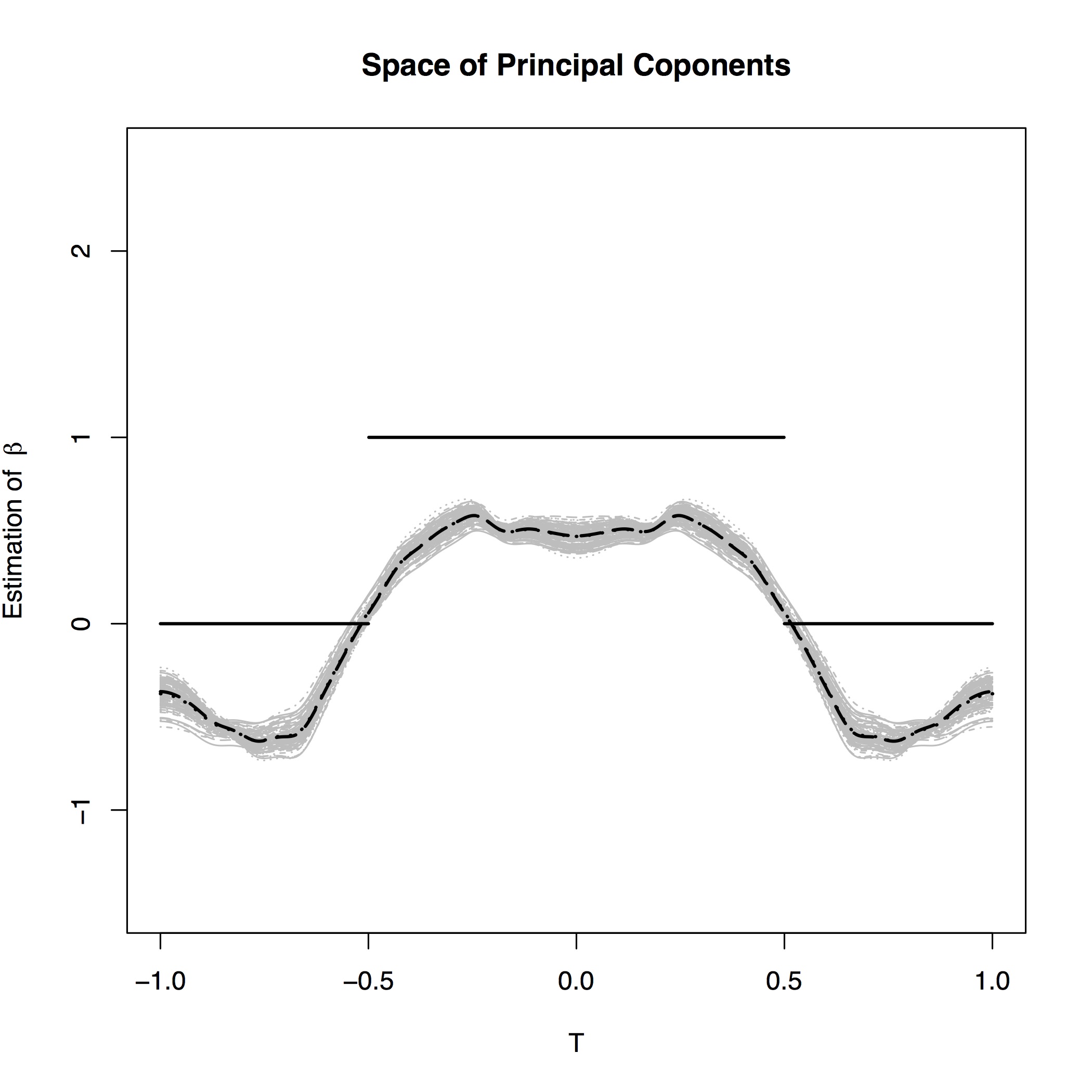}
    \caption{The black lines are the true $\beta(t)$, the dashed lines are the projection of $\beta(t)$ in the corresponding spaces,
    the dotted lines are the pointwise means of the estimated $\hat\beta(t)$.}
    \label{Simulazione_2}
\end{figure}

\textit{Case (c): $\beta\notin D$, and $\beta^F\neq0$.}
In this case, in general we have that $\gamma_n\neq0$ for $n\geq1$ and $\gamma\neq0$.
The asymptotic distance among $\widehat{\beta}^{D}_n$ and $\beta$ can be divided in three orthogonal terms:
$$\|\widehat{\beta}^{D}_n-\beta\|^2\ \stackrel{a.s.}{\rightarrow}\ \|\beta^F\|^2\ +\ \|\gamma\|^2\ +\ \|\left(\pi^{-1}\left(\beta^{E}+\gamma\right)-(\beta^E+\gamma)\right)-\beta^{S^{\perp}}\|^2,$$
where $\beta^F\in F$, $\gamma\in E$ and
$$\left(\pi^{-1}\left(\beta^{E}+\gamma\right)-(\beta^E+\gamma)\right)-\beta^{S^{\perp}}\ \in\ S^{\perp}.$$
In Figure~\ref{Simulazione_2}-left panel, we report $100$ independent simulations in which an estimate of $\widehat{\beta}^{D}_n$ is computed for large $n$
(see details in simulation setting in Appendix~\ref{appendix_simulation}).
Since in case (c) we are mainly interested in the estimation on $S$,
Figure~\ref{Simulazione_2} consider $D\equiv E$ and $\beta^{S^{\perp}}$, so that there is no bias on $S^{\perp}$.

The bias on the sub-space $F$ is always present in this situation, and it is simply due to the fact that
$\widehat{\beta}^{D}_n\in D$, which is included in $E\oplus S^{\perp}$, while $\beta\notin E\oplus S^{\perp}$ when $\beta^F\neq 0$.
Naturally, this bias also influences the statistical analysis on the outcome $y$, since the contribution of
$\langle X,\beta^F\rangle$ to $y$ is not taken into account.

When the aim of the analysis is to reconstruct only the component of $\beta$ on a particular sub-space,
given by $D$ and the functions orthogonal to the data, i.e. $E\oplus S^{\perp}$,
the bias on $F$ is not of interest.
However, the analysis on the estimation mainly focus on the bias on $E$: $\gamma(t)$.
This function indicates the asymptotic bias among $\widehat{\beta}^{E}_n$ and $\beta^{E}$.
%
In Figure~\ref{Simulazione_2}-left panel, $\gamma(t)$ is represented by the difference among $\beta$ (solid line) and
the pointwise mean of the estimated  $\hat\beta(t)$ (dotted line).
As mentioned in Subsection~\ref{section_estimation_betaE},
this bias is due to the fact that the part of the process $X$ on $E$ can be correlated to the part of $X$ along the component $\beta^{F}\in F$.
We may say that the bias $\gamma_n(t)$ puts in the estimate given by $\widehat{\beta}^{E}_n$ the additional information related to the contribution of $\langle X,\beta^{F}\rangle$ in computing $y_i$. Then, even if $\beta^{E}$ is the closest element of $E$ to $\beta$, the function of $E$ that better reconstructs $y_i$ from data is $(\beta^{E}+\gamma)$, since the contribution $\langle X,\beta^{F}\rangle$ is not observable.
From Remark~\ref{rem:bias_delta}, note that if $\bm{E}[X]=0$ and $E$ is composed by $d$ eigenfunctions of the covariance structure of $X$ (Karunen-Loeve basis),
then $\gamma(t)=0$.
In fact, in this case there is no information of $F$ contained in $E$, and then $\beta^{E}$ is also the function that better constructs $y_i$ from data in $E$.
In fact, in Figure~\ref{Simulazione_2}-right panel, where $D$ is the sub-space generated by the firsts principal components of the data,
the true $\beta$ (solid line) and
the pointwise mean of the estimated $\hat\beta(t)$ (dotted line) coincides ($\gamma(t)=0$).

\subsection{A bias-variance trade off in the estimation in finite sub-spaces}   \label{section_bias_variance}

In this subsection, we highlight an interesting bias-variance trade off concerning the choice of the sub-space
where the least square estimator is computed.
Before introducing this trade off, let us discuss the covariance structure of the estimator $\widehat{\beta}^{D}_n$,
which we define as $(\bm{\varphi^D(s)})^T\cdot\bm{Cov}(\bm{\widehat{\beta}^{D}_n})\cdot\bm{\varphi^D(t)}$,
since $\widehat{\beta}^{D}_n(t)=(\bm{\widehat{\beta}^{D}_n})^T\cdot\bm{\varphi^D(t)}$.
We now use the projection matrix $P$ such that,
$\bm{\widehat{\beta}^{D}_n}=P^{-1}(\bm{\widehat{\beta}^{E}_n})$, that is computed in Appendix~\ref{appendix_relation_D_E}.
Through this operator, we can express the relation among $\bm{Cov}(\bm{\widehat{\beta}^{D}_n})$ and $\bm{Cov}(\bm{\widehat{\beta}^{E}_n})$ as follows
$$\bm{Cov}(\bm{\widehat{\beta}^{D}_n})\ =\ P^{-1}\bm{Cov}(\bm{\widehat{\beta}^{E}_n})(P^{-1})^T.$$
From~\eqref{eq:P_E} in Appendix~\ref{appendix_relation_D_E}, we obtain
$$\bm{Cov}(\bm{\widehat{\beta}^{D}_n})\ =\ V_DD_D^{-1/2}\bm{Cov}(\bm{\widehat{\beta}^{E}_n})D_D^{-1/2}V_D^T,$$
where $D_D$ and $V_D$ represent the eigen-structure of $P^TP$, i.e. $P^TPV_D=V_DD_D$.
Denote with $\nu_1^{D},..,\nu_d^{D}$ and $\nu_1^{E},..,\nu_d^{E}$ the eigenvalues of $\bm{Cov}(\bm{\widehat{\beta}^{D}_n})$ and $\bm{Cov}(\bm{\widehat{\beta}^{E}_n})$, respectively. Then, we can observe that
\begin{itemize}
\item[(i)] since $P$ is a projection matrix, all the eigenvalues of $P^TP$ are less than one. Hence,
all the elements in $D_D^{-1/2}$ are greater than one. So, the variance of the retro-projection due to $P^{-1}$ is non decreasing in any direction, i.e. $\nu_k^{D}\geq \nu_k^{E}$ for any $k=1,..,d$;
\item[(ii)] if $D\subseteq S$, all the eigenvalues are equal to one and the total variance is the same, i.e. $\nu_k^{D}= \nu_k^{E}$ for any $k=1,..,d$;
\item[(iii)] if all the eigenvalues of $D_D$ are greater than a value $\epsilon_D>0$,
 we can uniformly control the variance of $\bm{\widehat{\beta}^{D}_n}$, i.e. $\nu_k^{D}\leq 1/\epsilon_D\cdot\nu_k^{E}$ for any $k=1,..,d$.
\end{itemize}
From these properties we can distinguish two interesting cases of bias-variance trade-off related to the choice of the sub-space $D$:
\begin{itemize}
\item[(1)] Consider all the possible identifiable sub-spaces $D$ with the same projection $E_0$ on $S$, i.e.
$D\subseteq E_0\oplus S^{\perp}$ and $D\cap S^{\perp}=0$.
From (i) and (ii) we have that the variance of $\widehat{\beta}^{D}_n$ is minimized by choosing $D\equiv E_0$, that is $D\subseteq S$.
However, to reduce the bias on $S^{\perp}$, some a priori information on $\beta$ may suggest another choice of $D$.
For instance, consider the case $\beta^F=0$ and $\beta^{S^{\perp}}\neq0$: if we choose $D\equiv E_0$ the variance of $\widehat{\beta}^{D}_n$ is minimized but we have a bias on $S^{\perp}$ (case (b)), while if we choose $D$ such that $\beta\in D$, the estimator $\widehat{\beta}^{D}_n$ has no bias but the variance may be very high.
Figure~\ref{Simulazione_1} describes this situation:
when $D\equiv E$ (Figure~\ref{Simulazione_1}-first panel),
the variance of the estimates is low but the pointwise mean of the estimated $\widehat{\beta}^D_n$ does not
target $\beta$; when $D=D_{\pi/6}$ (Figure~\ref{Simulazione_1}-second panel), the variance of the estimates increases and the bias decreases;
when $D=D_{\pi/3}$ (Figure~\ref{Simulazione_1}-third panel), the variance of the estimates is high but there is no bias.

When the sample size $n$ or the number of discrete observations $p$ are not too large,
we may prefer a small variance even if the estimator is biased.
Naturally, when we have no previous information on $\beta$, there is no chance to control the bias and the smartest choice
is to minimize the variance by choosing the closest $D$ to the space of the data $S$.

\item[(2)] Consider all the possible sub-spaces $E\subseteq S$, generated by the projection of $D$ on $S$.
To ease notation, take $D\subset S$ (i.e., $D\equiv E$).
It is well known that the variance of the estimator $\widehat{\beta}^{E}_n$ is smaller when the variance of the data is higher.
Then, the variance of $\widehat{\beta}^{D}_n$ is minimized when $D$ coincides with the space generated by the first principal components of $X$ on $S$.
However, to reduce the bias on $F$, some a priori information on $\beta$ may suggest a different choice of $D$.
For instance, taking $\beta^{S^{\perp}}=0$, when $D$ is equal to the space generated by the first principal components (PCs),
the variance of $\widehat{\beta}^{D}_n$ is minimized but we have a bias on $F$, since in general $\beta^F\neq0$ (case (c)); nevertheless, when $D$ is such that $\beta\in D$, the estimator $\widehat{\beta}^{D}_n$ has no bias (case (a)) but the variance may be very high.
Figure~\ref{Simulazione_3} describes this situation:
in Figure~\ref{Simulazione_3}-left panel, we have a space $D$ that includes $\beta$,
and so the estimates $\widehat{\beta}^D_n$ of $\beta$ are
unbiased but they show a large variance;
in Figure~\ref{Simulazione_3}-second panel, the space of the first PCs does not includes $\beta$,
and so the pointwise mean of the estimated $\widehat{\beta}^D_n$ does not target $\beta$,
but the variance is low.
\begin{figure}
\centering
  \includegraphics[width=0.4\textwidth]{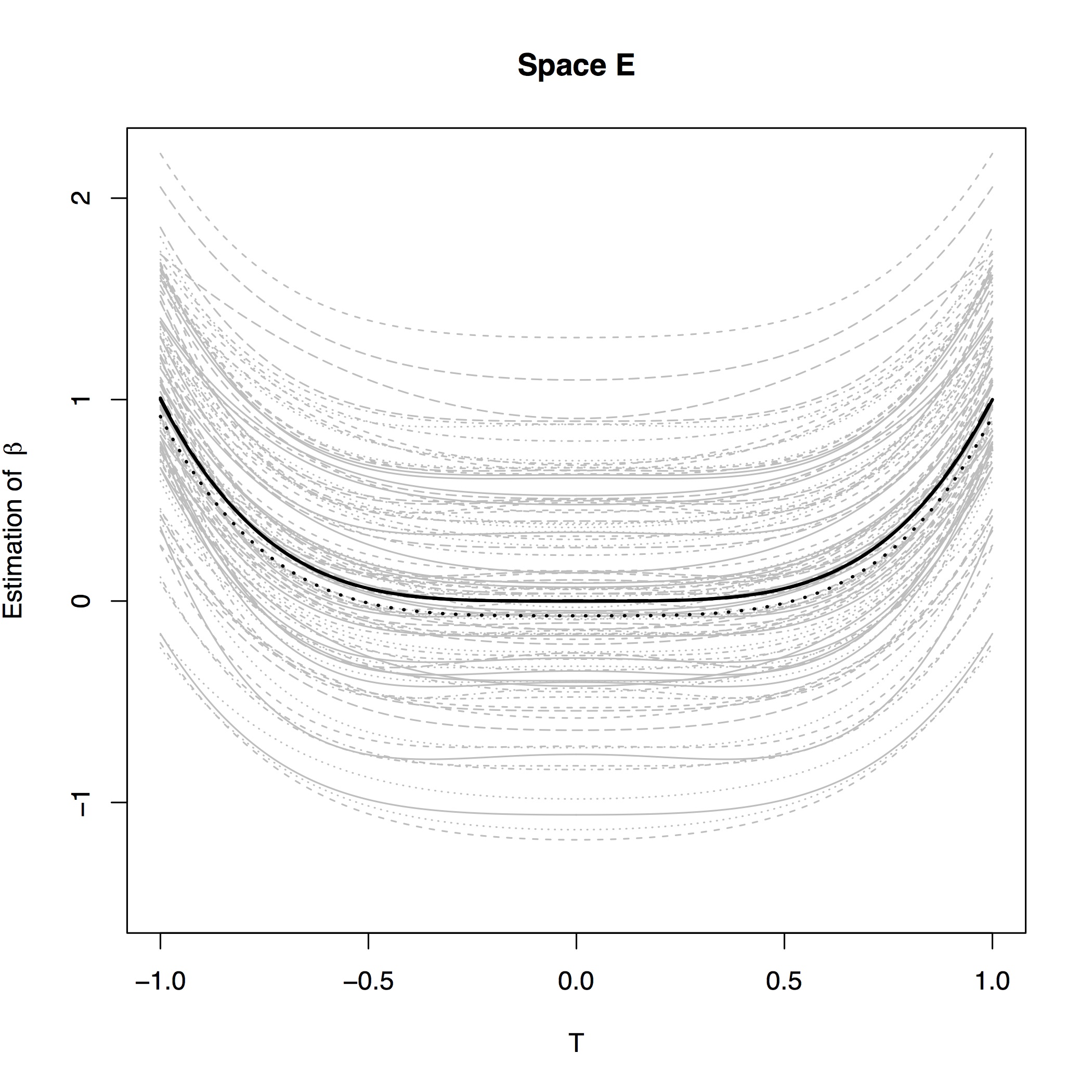}
  \includegraphics[width=0.4\textwidth]{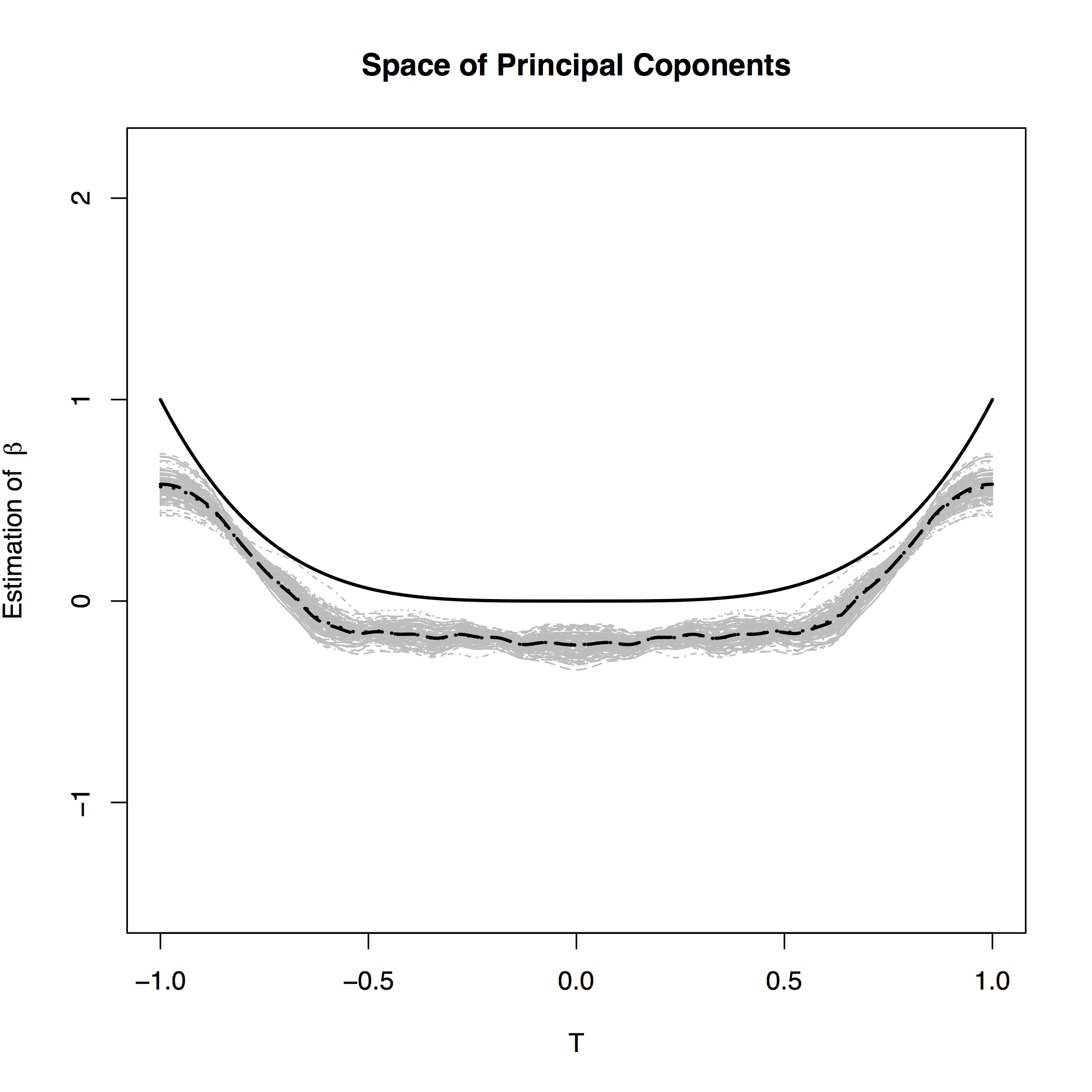}
   \caption{The black line is $\beta(t)$, the dashed line is the projection of $\beta(t)$ in the corresponding space, the dotted line is the pointwise mean of the estimated $\hat\beta(t)$. }
    \label{Simulazione_3}
\end{figure}

When the sample size $n$ or the number of discrete observations $p$ are not too large,
we may prefer a small variance even if the estimator is biased.
Naturally, when we have no a priori information on $\beta$, there is no chance to control the bias and the smartest choice
is to minimize the variance by choosing the closest $E$ to the space generated by the first PCs.
\end{itemize}

\section{Estimation in large dimensional sub-spaces}   \label{section_D_inf}

In this section, we discuss the behavior of the estimator $\widehat{\beta}_n^{D}$ obtained in~\eqref{eq:minimum_problem_D}
when the dimension of $D$ is arbitrarily large.
In other words, we want to investigate how to compute a well-defined estimator for $\beta$ in an infinite dimensional sub-space $D$.
To deal with this case, we express $D$ as the closure of a countable union of finite sub-spaces $\{D_d,d\geq1\}$, i.e.
$$D\ :=\ \overline{\bigcup_{d\geq1}\{D_d\}},$$
where $D_d\subset D_{d+1}$ for any $d\geq1$.
We denote with $\{\varphi_d^D,d\geq1\}$ the orthonormal basis of $D$, such that $\{\varphi_k^D,k=1,..,d\}$ is an orthonormal basis of $D_d$,
for any $d\geq1$. Note that $\dim(D_d)=d$.
A basic idea to construct an estimation procedure of $\beta$ in $D$ is to consider the estimators
$\{\widehat{\beta}_n^{D_d};d\geq1\}$ computed in the finite dimensional spaces $\{D_d;d\geq1\}$ by~\eqref{eq:minimum_problem_D},
and investigate their asymptotic behavior for large $d$.
In fact, from~\eqref{eq:convergence_beta_D} we have that $\{\lim_{n\rightarrow\infty}\widehat{\beta}_n^{D_d}\}$ exists finite for any fixed $d\geq1$;
however, $\widehat{\beta}_n^{D_d}$ can be considered a proper estimator in $D$ for $\beta$ only if the sequence of the limits $\left\{\lim_{n\rightarrow\infty}\widehat{\beta}_n^{D_d};d\geq1\right\}$ is convergent when $d\rightarrow\infty$.

Here we consider sub-spaces $D_d$ with arbitrarily large dimension,
it is worth making an important consideration on the estimator $\widehat{\beta}_n^{D_d}$.
As mentioned in Subsection~\ref{section_least square}, a least square estimator for $\beta$ in a finite identifiable sub-space is well-defined only
if both the sample size $n$ and the number of observations per curve $p$ are greater than the dimension of the sub-space itself.
Then, $\widehat{\beta}_n^{D_d}$ can be computed only if $\min\{n;p\}\geq d$;
moreover, whenever we let $d$ increase to infinity,
we are implicitly requiring that both $n$ and $p$ must diverge with a rate depending on $d$.
Therefore, in all the situations in which $n$ or $p$ can't increases arbitrarily,
the results presented in this section do not hold.

In the following, we consider a framework analogous to the one presented in Section~\ref{section_estimation_subspace}:
for each $d\geq1$, let $E_d$ be the sub-space obtained by the projection of $D_d$ on $S$, i.e.
$$E_d\ :=\ \left\{\ x\in S\ :\ \exists y\in D_d, x=\sum_{k=1}^{\dim(S)}\langle y,\varphi^S_k\rangle\varphi^S_k\ \right\},$$
and let define
$$E\ :=\ \overline{\bigcup_{d\geq1}\{E_d\}},\ \ \ \ F_d:=S\cap E_d^{\perp},\ \ \ \ \ F:=S\cap E^{\perp}.$$
So, we have that $D\subset E \oplus S^{\perp}$ and
$L^2(T)= E \oplus F \oplus S^{\perp}$ and any $\beta\in L^2(T)$ has the following orthogonal decomposition: $\beta= \beta^E + \beta^F + \beta^{S^{\perp}}$.

\subsection{Estimation instability in large dimensional sub-spaces}   \label{subsection_controesempio}

In this subsection, we show that the limit of the sequence $\left\{\lim_{n\rightarrow\infty}\widehat{\beta}_n^{D_d};d\geq1\right\}$ may not exist,
even when $\beta\in D$.
To do this, we discuss an example where the sequence
$\left\{\lim_{n\rightarrow\infty}\widehat{\beta}_n^{D_d};d\geq1\right\}$
does not converge when $d\rightarrow\infty$.
Since from~\eqref{eq:convergence_beta_D} $\widehat{\beta}_n^{D_d}\rightarrow_n\beta^{D_d}+\gamma^d$ a.s. for any $d\geq1$,
where $\gamma^d$ is the asymptotic bias on $E_d$ and since $\beta^{D_d}\rightarrow_d\beta^{D}$, then
it is sufficient to show that $\{\|\gamma^d\|;d\geq1\}$ is not bounded as $d$ increases.\\

Consider a process $X$, with $\bm{E}[X]=0$, defined on an infinite dimensional sub-space $D$,
with Karhunen-Loeve (K-L) basis $\{\psi_k;k\geq1\}$ and corresponding eigenvalues $\{\lambda_k;k\geq1\}$. The sequence $\{\lambda_k;k\geq1\}$ is decreasing in $k$ (i.e. $\lambda_{max} \equiv \lambda_1 \geq \lambda_2 \geq \cdots$).
Let $\{\varphi_k;k=1,..,d\}$ be a basis for $E_d$, defined as follows:
\begin{equation}\label{eq:phi_controesempio}
\begin{aligned}
\varphi_k & =
\begin{cases}
\cos(\theta_k)\psi_k+\sin(\theta_k)\psi_{k+1} & \text{if }k \text{ odd};
\\
\sin(\theta_{k-1})\psi_{k-1}-\cos(\theta_{k-1})\psi_{k} & \text{if }k \text{ even},
\end{cases}
\end{aligned}
\end{equation}
where the sequence $\{\theta_k;k\geq1\}$ will be appropriately determined more ahead.
Using the representation $\bm{\gamma^d}=V^E\bm{\delta^d}$ presented in Remark~\ref{rem:bias_delta},
it is sufficient to show that $\{\|\delta^d\|;d\geq1\}$ is not bounded as $d$ increases.
To this aim, note that the K-L basis of $X$ projected on $E_d$ is $\{\psi_k;k=1,..,d-1\}\cup\{\varphi_d\}$ for $d$ odd, and
$\{\psi_k;k=1,..,d\}$ for $d$ even.
By Remark~\ref{rem:bias_delta} it is easy to see that $\delta^d_k=0$ for any $k<d$,
while $\delta^d_d=0$ when $d$ is even and
$$\delta^d_d\ =\ \frac{\bm{Cov}\left(\langle X,\varphi_d\rangle,\langle X, \beta^{F_d}\rangle\right)}
{\bm{Var}\left(\langle X,\varphi_d\rangle\right)},$$
when $d$ is odd.
Hence, $\|\delta^d\|=0$ for $d$ even, while $\|\delta^d\|=|\delta^d_d|$ for $d$ odd.
This last term is not zero because of the correlation among the projection of $X$ on $\varphi_d$ (included in $E_d$) and the projection of $X$ on $\varphi_{d+1}$ (included in $F_d$).
By writing $\beta=\sum_{k\geq1}\beta_k\psi_k$, and $X$ projected on $E_d$ as $\sum_{k=1}^{d}Z_k\sqrt{\lambda_k}\psi_k$, we obtain
$$\langle X,\varphi_d\rangle\ =\ Z_d\sqrt{\lambda_d}\cos(\theta_d)+Z_{d+1}\sqrt{\lambda_{d+1}}\sin(\theta_d),$$
\[\begin{aligned}
\langle X, \beta^{F_d}\rangle\ &&=&\ \sum_{k=d+2}^{\infty}Z_k\beta_k\sqrt{\lambda_k}\\
&&+&\ \tilde{\beta}_d
\left(Z_d\sqrt{\lambda_d}\sin(\theta_d)-Z_{d+1}\sqrt{\lambda_{d+1}}\cos(\theta_d)\right),
\end{aligned}\]
where $\tilde{\beta}_d=\left|\beta_d\sin(\theta_d)-\beta_{d+1}\cos(\theta_d)\right|$.
Then, from some easy calculations we have that
$$|\delta_d^d|\ =\ \left(\frac{|\cos(\theta_d)\sin(\theta_d)|\mu_d}
{\mu_d\cos^2(\theta_d)+1}\right)\cdot\tilde{\beta}_d,$$
where $\mu_d=\lambda_d/\lambda_{d+1} - 1$.
Now, consider any sequence $\{\epsilon_d;d\geq1\}$ such that $|\beta_d|/\epsilon_d\rightarrow\infty$, and take $\theta_d=\pi/2-\epsilon_d$ and $\lambda_{d+1}=\lambda_d/(1+\exp(\epsilon_d^{-1}))$, so that for large odd $d$
$$\|\delta^d\|\ \simeq\ \frac{|\beta_d|}{\epsilon_d(1+\exp(-\epsilon_d^{-1}))}\ \rightarrow_d\ \infty.$$
This concludes the example where the sequence $\{\lim_{n\rightarrow\infty}\widehat{\beta}_n^{D_d};d\geq1\}$ does not converge for $d\rightarrow\infty$.

\subsection{Principal components for estimation in large dimensional sub-spaces}   \label{subsection_D_inf_abc}

In Subsection~\ref{subsection_controesempio} we have shown that the limit of the sequence
$\left\{\lim_{n\rightarrow\infty}\widehat{\beta}_n^{D_d};d\geq1\right\}$
does not exist in general, not even when $\beta\in D$.
We discuss how to introduce an alternative least square estimator well defined in the case of $\beta\in D$.
We will denote this estimator as $\widetilde{\beta}_n^{d,k}$,
where $n\geq1$ is the sample size and $d\geq k\geq1$ are two integer parameters associated to the dimension of the sub-space.
In this subsection, we show that, when $\beta\in D$, there exists a sequence $k_d\rightarrow\infty$ such that
the sequence $\left\{\lim_{n\rightarrow\infty}\widetilde{\beta}_n^{d,k_d};d\geq1\right\}$ converges when $d\rightarrow\infty$.
This will let us consider $\widetilde{\beta}_n^{d,k_d}$ as a proper estimator for $\beta$ when $d$ is large.
To obtain this result, we need to assume the following conditions
\begin{itemize}
\item[(i)] $\beta\in E \oplus S^{\perp}$, that means $\beta^F=0$;
\item[(ii)] $D\subseteq S$, that means $D\equiv E$ and $D_d\equiv E_d$ for any $d\geq1$.
\end{itemize}
It is worth highlighting that these conditions are not restrictive and in literature they are always assumed to be true.
In fact, most of the existent works consider the limiting space $D$ equal to the space that generates the data, i.e. $D\equiv S$,
which implies both conditions (i) and (ii).\\

Let $\{\psi^d_i;i=1,..,d\}$ be the K-L basis of $X$ projected on the sub-space $D_d$, for any $d\geq1$ and
recall that $\widehat{\beta}_n^{D_d}$ is the least square estimator computed on $D_d$ from~\eqref{eq:minimum_problem_D};
then, we define $\widetilde{\beta}_n^{d,k}$ as the projection of $\widehat{\beta}_n^{D_d}$ on the sub-space
generated by the first $k$ functions of the K-L expansion in $D_d$, i.e.
\begin{equation}\label{eq:beta_tilde}
\widetilde{\beta}_n^{d,k}(t)\ :=\ \sum_{i=1}^k\langle\widehat{\beta}_n^{D_d},\psi^d_i\rangle\psi^d_i(t).
\end{equation}
Analogously, we define $\beta^{d,k}$ and $\gamma^{d,k}$ as the projections of $\beta$ and $\gamma^d$, respectively,
on the sub-space generated by the first $k$ eigenfunctions of $X$ in $D_d$, i.e.
\begin{equation}\label{eq:beta_gamma_dk}
\beta^{d,k}(t):=\sum_{i=1}^k\langle\beta,\psi^d_i\rangle\psi^d_i(t),\ \ \ \ \ \
\gamma^{d,k}(t):=\sum_{i=1}^k\langle\gamma^{d},\psi^d_i\rangle\psi^d_i(t).
\end{equation}
Since from~\eqref{eq:convergence_beta_D} we have that $\widehat{\beta}_n^{D_d}\rightarrow_n \beta^{D_d}\ +\ \gamma^{d}$ a.s.,
we can project all the terms on the sub-space generated by $\{\psi_1^d,..,\psi_k^d\}$, obtaining
$$\widetilde{\beta}_n^{d,k}\ \stackrel{a.s.}{\rightarrow}\ \beta^{d,k}\ +\ \gamma^{d,k}.$$

It is trivial to show that $\beta^{d,k}\rightarrow \beta^D$ when $d$ and $k$ increase to infinity,
then, our aim is to show that there exists a sequence $k_d\rightarrow\infty$ such that
\begin{equation}\label{eq:gamma_kd_tends_to_zero}
\|\gamma^{d,k_d}\|\ \rightarrow\ 0.
\end{equation}
To do that, fix $k\leq d$ and consider the coefficients of $\gamma^{d,k}$ with respect the basis $\{\psi_1^d,..,\psi_k^d\}$, i.e.
$\delta_i^d=\langle\gamma^d,\psi_i^d\rangle$ for $i=1,..,k$, where from Remark~\ref{rem:bias_delta}
$$\delta_i^d\ =\ \left(\bm{Var}\left[\langle X,\psi_i^d\rangle\right]\right)^{-1}\cdot
\bm{Cov}\left[\langle X,\psi_i^d\rangle,\langle X, \beta^{F}\rangle\right].$$
Then, defining $\lambda_i^d:=\bm{Var}\left[\langle X,\psi_i^d\rangle\right]$ and applying Cauchy-Schwartz inequality,
we obtain that, for any $i=1,..,k$,
$$(\delta_i^d)^2\ \leq\ (\lambda_i^d)^{-1}\bm{Var}\left[\langle X, \beta^{F}\rangle\right]\ \leq\
\left(\frac{\lambda_{max}}{\lambda_i^d}\right)\|\beta^{F_d}\|^2.$$
From~\eqref{eq:Cauchy_for_paper} we have that $\lambda_i^d$ is increasing in $d$, so that $\lambda_i^d\geq\lambda_i^k$.
Therefore, for any $i=1,..,k$, we have that
$$
(\delta_i^d)^2\ \leq\ \left(\frac{\lambda_{max}}{\lambda_i^k}\right)\|\beta^{F_d}\|^2,
$$
and hence
$$\|\gamma^{d,k}\|^2\ =\ \sum_{i=1}^k\left(\delta_i^d\right)^2\ \leq\
\sum_{i=1}^k\left(\frac{\lambda_{max}}{\lambda_i^i}\right)\|\beta^{F_d}\|^2.$$
Moreover, from~\eqref{eq:Cauchy_for_paper} we know that $\lambda_k^k\leq\lambda_i^i$ for any $i\leq k$,
and calling $C_k=k(\lambda_{max}/\lambda_k^k)$, we obtain
\begin{equation}\label{eq:bias_PC}
\|\gamma^{d,k}\|^2\ \leq\
k\left(\frac{\lambda_{max}}{\lambda_k^k}\right)\|\beta^{F_d}\|^2\ =\ C_k\|\beta^{F_d}\|^2,
\end{equation}
for any fixed $k\geq1$.
Since $\|\beta^{F_d}\|\rightarrow_d0$ because $\beta^F=0$,
we can take a sequence $k_d\rightarrow_d \infty$ such that
$C_{k_d}\|\beta^{F_d}\|^2\rightarrow0$, so that from~\eqref{eq:bias_PC} we get~\eqref{eq:gamma_kd_tends_to_zero}.
As a consequence, the sequence $\left\{\lim_{n\rightarrow\infty}\widetilde{\beta}_n^{d,k_d};d\geq1\right\}$ converges when $d\rightarrow\infty$,
which let us consider $\widetilde{\beta}_n^{d,k}$ as a proper estimator of $\beta\in D$ for large $d$.\\

Finally, we can write the consistency result for the estimator $\widetilde{\beta}_n^{d,k}$, by letting $k$ and $d$ depending on the sample size $n$:
under assumptions (i) and (ii), there exists a sequence $\{d_n;n\geq1\}$ such that
\begin{equation}\label{eq:consistency_D_inf}
\widehat{\beta}^{d,k}_n\ \stackrel{a.s.}{\rightarrow}_n\ \beta^D,
\end{equation}
where $d=d_n$ and $k=k_{d_n}$ for any $n\geq1$.
Result~\eqref{eq:consistency_D_inf} can be written as follows
$$\|\widetilde{\beta}^{d,k}_n-\beta\|\ \stackrel{a.s.}{\rightarrow}_n\ \|\beta^{S^{\perp}}\|,$$
which implies that, when $\beta\in D$,
$$\|\widetilde{\beta}^{d,k}_n-\beta\|\ \stackrel{a.s.}{\rightarrow}_n\ 0.$$

\begin{rem}\label{rem:condition_essential}
Assumption (i) is essential to consider $\widehat{\beta}^{d,k}_n$ as a proper estimator of $\beta$.
To see this, consider the following example, where (i) fails, i.e. $\beta^F\neq0$, and there is no sequence $\{k_d;d\geq1\}$
such that $\|\gamma^{d,k_d}\|$ is convergent.
In particular, let $\beta=c\cdot\phi$, with $\phi\in F$ and $\|\phi\|=1$.
Then, take a process $X$ defined as follows
$$X\ =\ \sum_{d=1}^{\infty}Z_d\sqrt{\lambda_d}\varphi_d\ +\ \left(\sum_{j=1}^{\infty}Z_j\sqrt{\lambda_j}\right)\phi$$
where $\{\varphi_d;d\geq1\}$ is an orthonormal basis of $D$ and $\{Z_k;k\geq1\}$ are i.i.d. r.v. with zero mean and unit variance.
Then, define the sequence $\{D_d;d\geq1\}$ as $D_d=\textit{span}\{\varphi_1,..,\varphi_d\}$.
Hence, for any $k=1,..,d$ we have that $\gamma_k^d=1$, which implies $\|\gamma^{k_d,d}\|=\sqrt{k_d}\rightarrow_d\infty$ for any
divergent sequence $\{k_d;d\geq1\}$.
\end{rem}

\appendix
\section{Formal characterization of the sub-space $E$}   \label{appendix_relation_D_E}

This section focuses on computing explicitly the following quantities introduced in the Section \ref{section_estimation_subspace}:
\begin{itemize}
\item[(1)] the orthonormal basis of $E$: $\left\{\varphi_k^{E};k=1,..,d\right\}$;
\item[(2)] the multivariate projection matrix $P:\mathbb{R}^d\rightarrow \mathbb{R}^d$ that transforms the basis coefficients of elements in $D$
in the basis coefficients of elements in $E$.
\item[(3)] the functional projection operator $\pi:D\rightarrow E\subseteq S$ of $D$ on $S$
\end{itemize}
Let us consider point (1). First, project the basis of $D$ ($\left\{\varphi_k^{D};k=1,..,d\right\}$) on $S$, so obtaining a $\dim(S)\times d$-matrix $A$, where $[A]_{ij}=\langle\varphi_i^S,\varphi_j^D\rangle$. Note that $A$ may have infinite rows if $\dim(S)=\infty$.
Then, the basis of $D$ projected on $S$ generates $d$ linear independent functions given by $A^T\bm{\varphi^S(t)}$, that is a basis for $E$.
It is easy to show that $A^T\bm{\varphi^S(t)}$ are linear independent since $\varphi_1^{D},..,\varphi_d^{D}$ are, and $D\cap S^{\perp}=0$.
To make $A^T\bm{\varphi^S(t)}$ be an orthonormal basis for $E$ we do some calculations, obtaining:
\begin{equation}\label{eq:basis_E}
\bm{\varphi^{E}(t)}=V_SD_D^{-1/2}V_D^TA^T\bm{\varphi^S(t)},
\end{equation}
where $D_D$ and $V_D$ represent the eigen-structure of $A^TA$ ($A^TAV_D=V_DD_D$) and $V_S$ is an arbitrary $d\times d$-orthonormal matrix
that allows the basis of $E$ to be changed; without loss of generality, we can consider $V_S=I_d$.
Note that, except for $V_S$, the basis $\bm{\varphi^{E}(t)}$ is independent of the choice of the basis $\bm{\varphi^{D}(t)}$ and $\bm{\varphi^{S}(t)}$.
It is worth saying that the eigenvalues in $D_D$ are all strictly positive since $A^TA$ has full rank,
since $\varphi_1^{E},..,\varphi_d^{E}$ are linear independent.
Moreover, the eigenvalues in $D_D$ are all less or equal to one since $A$ is a projection operator.

Now, consider point (2).
From~\eqref{eq:basis_E} the projection matrix $P$ from $D$ to $E$ can be defined as
$$P\ :=\ \langle \bm{\varphi^E(t)},\left(\bm{\varphi^D(t)}\right)^T\rangle\ =\
V_SD_D^{-1/2}V_D^TA^T\langle \bm{\varphi^S(t)},\left(\bm{\varphi^D(t)}\right)^T\rangle,$$
since $\langle \bm{\varphi^S(t)},\left(\bm{\varphi^D(t)}\right)^T\rangle=A$ and
$V_D^TA^TA=D_DV_D^T$, we obtain
\begin{equation}\label{eq:P_E}
P=V_SD_D^{1/2}V_D^T.
\end{equation}
Note that, using~\eqref{eq:P_E} we can rewrite~\eqref{eq:basis_E} as
$$\bm{\varphi^{E}(t)}=(P^{-1})^TA^T\bm{\varphi^S(t)}.$$
Then, from the vectorial estimate in $E$ given by~\eqref{eq:minimum_problem_E}, we can obtain the vectorial estimate in $D$ with $\bm{\widehat{\beta}^{D}_n}=P^{-1}(\bm{\widehat{\beta}^{E}_n})$, and finally compute the functional estimate
$\widehat{\beta}^{D}_n(t)=(\bm{\widehat{\beta}^{D}_n})^T\bm{\varphi^{D}(t)}$. This coincides with the solution of~\eqref{eq:minimum_problem_D}.

Finally, consider point (3).
Using the projection matrix $P$ we can define the functional operator $\pi$ as follows
$$\pi(g)=\left(P\langle g,\bm{\varphi^{D}(t)}\rangle\right)^T\bm{\varphi^{E}(t)},$$
for any $g\in D$. Then, using~\eqref{eq:P_E} we can easily obtain
\begin{equation}\label{eq:pi_E}
\pi(\cdot)=\left(\langle \cdot,\bm{\varphi^{D}(t)} \rangle\right)^T A^T\bm{\varphi^{S}(t)}.
\end{equation}
Note that $\pi$ is independent of any choice of basis of $S$, $D$ and $E$.
Using~\eqref{eq:pi_E}, once we get the vectorial estimate in $E$ from~\eqref{eq:minimum_problem_E}, we can immediately compute the functional estimate
$\widehat{\beta}^{E}_n(t)=(\bm{\widehat{\beta}^{E}_n})^T\bm{\varphi^{E}(t)}$, and then obtain the functional estimate in $D$, i.e.
$\widehat{\beta}^{D}_n=(\pi)^{-1}(\widehat{\beta}^{E}_n)$.

\section{Increasing information property}   \label{appendix_cauchy}

In this section, we discuss an interesting property concerning the behavior of the eigenvalues of the covariance matrix when its dimension increases.

Let $\{M^{(n)}=[m^{(n)}_{ij}],n\geq1\}$ be a sequence of symmetric
matrices such that,
for each $n\geq1$, $M^{(n)}$ is a $n\times n$ matrix with $m^{(n)}_{ij}=m^{(n-1)}_{ij}$ for any $i,j\leq n-1$.
In other words, $M^{(n-1)}$ is obtained by $M^{(n)}$ by deleting the lat row and column.
The eigenvalues are real, and are ordered according to the following general result proved by
Cauchy in~\cite[p.~187]{Cauchy.29}.
\begin{thm}\label{lem:Chauchy}
{\bf(\cite[p.~125]{Hawkins.77})}
On the nested sequence $(M^{(n)})_n$ of matrices given above,
denote with $\{\lambda^{n}_k;k=1..,n\}$ the sequences of the ordered eigenvalues of $M^{(n)}$.
Then, for any $n\geq 1$,
$$
\lambda^{n+1}_1\ \geq\ \lambda^{n}_1\ \geq\ \lambda^{n+1}_2\ \geq\ \lambda^{n}_2\ \geq\ \lambda^{n+1}_3\ \geq
\cdots \geq\ \lambda^{n}_{n}\ \geq\ \lambda^{n+1}_{n+1}.
$$
\end{thm}
A direct consequence of the previous theorem is
\begin{equation}\label{eq:Cauchy_for_paper}
\lambda_i^{k}\leq\lambda_i^{d}, \quad \lambda_i^{i} \leq \lambda_k^{k}, \qquad \forall i\leq k\leq d.
\end{equation}
This result is applied in Section~\ref{subsection_D_inf_abc},
where $M^{(n)}$ is represented the covariance matrix of the random vector $(\langle X,\varphi_1\rangle,..,\langle X,\varphi_n\rangle)$.
In this context, a direct interpretation of~\eqref{eq:Cauchy_for_paper} is that
the variance of $X$ projected into a subspace increases when further components are added.

\section{Simulation settings}   \label{appendix_simulation}
The settings of the simulation study presented in Section~\ref{section_estimation_subspace} are the following.
\noindent
\begin{itemize}
\item[(1)]
Data $x_i(t)$ and regression coefficient $\beta(t)$ belong to the Hilbert space $L^2(T)$ with $T= [-1,1]$ closed interval.
\item[(2)]
The finite dimensional sub-spaces we consider are:
$$E = \text{Span} \{1/\sqrt{2}, \sqrt{5/8}(3t^2-1), \sqrt{9/128} (35 t^4 - 30 t^2 +3)\}$$
and
$$D_{\theta} = \text{Span} \{\cos{(\theta)} 1/\sqrt{2} + \sin{(\theta)}\sqrt{3/2} t,  \sqrt{5/8}(3t^2-1), \sqrt{9/128} (35 t^4 - 30 t^2 +3)\},$$ with $\theta \in [0,2\pi]$.
\end{itemize}
\noindent
Observe that $E \equiv D_0$.

\noindent
For each $i = 1,...,n$ where $n$ is the sample size (in our examples $n=500$),
$$X_i(t) = \sum_{j \in J_i} \alpha_j \eta_j \theta^{X}_j(t), $$ where
$\{\theta_k^{X}(t)\} \equiv \{1/\sqrt{2}\} \bigcup \{\cos{(\pi k t)}, k = 1,...\}$, $\alpha_j$ are randomly sampled from a uniform distribution
$U \sim \text{Unif}_{[-10,10]}$, $\eta_1 = 0.01, \eta_j = 1/j$, $for j > 1$ and
$J_i$ is a subset of size $Z$ (with $Z$ Poisson random variable $Z \sim \mathcal{P}(\lambda)$) of the integer from $1$ to $2*Z$.
We set $\lambda = 10$.

Chosen a function $\beta(t) \in L^2(T)$  the scalar responses $y_1,...,y_n$ are generated as $y_i = \int_T \beta(t)X_i(t) dt + \epsilon_i$, where $\epsilon_i \sim \mathcal{N}(0,1)$.
We repeat the estimation procedure $M = 100$ times.

In Figure~\ref{Simulazione_1} the true $\beta(t)$ is $\beta(t) = t^2 + 2t + 1/3$,
in Figure~\ref{Simulazione_2}: the true $\beta(t)$ is $ \beta(t)= \ind_{[-0.5,0.5]}(t)$ and in Figure~\ref{Simulazione_3} the true $\beta(t)$ is $ \beta(t) = t^4$.

\section*{Acknowledgements}
The authors wish to thank Piercesare Secchi for stimulating and essential discussions about topics covered by this paper.

\end{document}